\documentclass[a4paper, 11pt]{article}
\title{On pseudo-null Iwasawa modules, restricted ramification, and Non-freeness Conjecture over a quartic CM-field}
\author{Takuya Yanagisawa\\
{\small Department of Mathematics, Graduate School of Science, Tokyo University of Science}\\
{\small \texttt{1125705@ed.tus.ac.jp}}}
\date{}

\usepackage{amsmath, amssymb}
\usepackage{amsthm}
\usepackage{mathtools}
\usepackage{stmaryrd}
\usepackage{float}
\usepackage{url}
\usepackage[labelsep=none]{caption}
\usepackage{enumitem}

\usepackage[top=3cm,bottom=3cm,left=2cm,right=2cm]{geometry}

\theoremstyle{plain}
\newtheorem{df}{Definition}[section]
\newtheorem{thm}[df]{Theorem}
\newtheorem{prop}[df]{Proposition}

\newtheorem{lem}[df]{Lemma}

\newtheorem{conjecture}[df]{Conjecture}

\theoremstyle{definition}
\newtheorem{rem}[df]{Remark}
\newtheorem{example}[df]{Example}

\usepackage{tikz}
\usetikzlibrary{arrows.meta}

\newcommand{\Z}{\mathbb{Z}}
\newcommand{\Q}{\mathbb{Q}}

\newcommand{\F}{\mathbb{F}}

\newcommand{\pe}{\mathfrak{p}}
\newcommand{\Pe}{\mathfrak{P}}
\newcommand{\qu}{\mathfrak{q}}

\newcommand{\calG}{\mathcal{G}}
\newcommand{\calL}{\mathcal{L}}

\DeclareMathOperator{\Ker}{Ker}
\DeclareMathOperator{\Gal}{Gal}
\DeclareMathOperator{\Ann}{Ann}
\DeclareMathOperator{\rank}{rank}
\DeclareMathOperator{\Homo}{H}
\DeclareMathOperator{\Tor}{Tor}

\DeclareMathOperator{\pd}{pd}
\DeclareMathOperator{\cd}{cd}
\DeclareMathOperator{\chara}{char}
\DeclareMathOperator{\length}{length}

\newcommand{\ab}{\mathrm{ab}}
\newcommand{\cy}{\mathrm{c}}

\newcommand{\tgen}[1]{\overline{\langle #1 \rangle}}
\newcommand{\br}[1]{\lbrack #1 \rbrack}
\newcommand{\bbr}[1]{\llbracket #1 \rrbracket}

\newcommand{\wtk}{\widetilde{k}}

\usepackage[backend=biber, style=numeric]{biblatex}
\begin{document}
\maketitle

\begin{abstract}
Let $k$ be a CM-field of degree greater than $2$, $p$ an odd prime number that splits completely in $k$, and $\Sigma$ a $p$-adic CM-type of $k$. In this setting, we show the non-triviality of the $\Sigma$-ramified Iwasawa module of the maximal multiple $\Z_p$-extension $\wtk/k$. As a consequence, under a weak form of Greenberg's Generalized Conjecture, we see that the unramified Iwasawa module and the $\Sigma$-ramified Iwasawa module of $\wtk$ never coincide. In the case where $k$ is a quartic CM-field, under a certain assumption, we show that the non-coincidence condition is equivalent to Greenberg's  Generalized Conjecture. We also study a relation of the non-coincidence condition with Non-freeness Conjecture, which is a conjecture studied in non-abelian Iwasawa theory. 
\end{abstract}

\section{Introduction}
Let $p$ be a prime number. For an algebraic extension $K/\Q$, we let $L(K)/K$ denote the maximal unramified abelian pro-$p$-extension and $X(K)\coloneqq\Gal(L(K)/K)$ its Galois group. If $K$ is a Galois extension over a certain subfield $k$, then $X(K)$ naturally has a module structure over the completed group algebra $\Z_p\bbr{\Gal(K/k)}$. In the case where $K/k$ is a $\Z_p^d$-extension for some $d\geq 1$, the module $X(K)$ is often called the unramified Iwasawa module of $K/k$. 

We also consider Iwasawa modules with restricted ramification. Let $S_p(k)$ be the set of $p$-adic primes of $k$ and $\Sigma$ an element or a subset of $S_p(k)$. In this setting, for any algebraic extension $K/k$, we let $M_\Sigma(K)/K$ denote the maximal abelian pro-$p$-extension that is unramified outside primes lying above $\Sigma$ (we refer to this property as $\Sigma$-ramified for short). We put $X_\Sigma(K)\coloneqq\Gal(M_\Sigma(K)/K)$. If $K/k$ is a Galois extension, we see that $X_\Sigma(K)$ has a module structure over $\Z_p\bbr{\Gal(K/k)}$, and if $K/k$ is a $\Z_p^d$-extension for some $d\geq 1$, the module is called the $\Sigma$-ramified Iwasawa module of $K/k$. Note that $L(K)$ (resp. $X(K)$) coincides with $M_\emptyset(K)$ (resp. $X_\emptyset(K)$). 

In the present article, we study the relations between four conjectures in classical Iwasawa theory, including a new one.

\subsection{Non-triviality of Iwasawa modules}
For a finite extension $k/\Q$, we let $\wtk/k$ denote the maximal multiple $\Z_p$-extension, that is, the composite of all $\Z_p$-extensions over $k$. Concerning the extension $\wtk/k$, the following conjecture has been studied. 
\begin{conjecture}[Greenberg's Generalized Conjecture (GGC) {\cite[Conjecture~(3.5)]{Greenberg2001}}]\label{conjecture GGC}
Let $k/\Q$ be a finite extension and $p$ a prime number. Then $X(\wtk)$ is pseudo-null as a $\Z_p\bbr{\Gal(\wtk/k)}$-module. 
\end{conjecture}
A useful overview of earlier works on GGC can be found in the introduction of \cite{Murakami2023}. 

Concerning GGC, an extreme case is also studied, namely, the case where $X(\wtk)$ is trivial. For example, it is known that the triviality of $X(\wtk)$ places a strong restriction on the degree of the base field as follows. 

\begin{prop}[{\cite[Lemma~3.9]{Fujii2022}}]\label{prop Fujii non-trivial}
Let $k/\Q$ be a totally imaginary finite extension, $p$ a prime number that splits completely in $k$. If $[k:\Q]\geq 8$, then $X(\wtk)$ is not trivial. 
\end{prop}
\begin{rem}\label{rem Leopoldt}
The validity of $p$-adic Leopoldt's Conjecture is assumed in \cite[Lemma~3.9]{Fujii2022}, but it is not necessary. See \S\ref{subsection unramified} for details. 
\end{rem}
In the case where $k$ is a CM-field, a more detailed result is shown by Jaulent in \cite[Theorem~10]{Jaulent2024}. See recent preprints \cite{Okano2026pre} and \cite{MO2026pre}, as well. 
We show the following variant of Proposition~\ref{prop Fujii non-trivial}. 
\begin{prop}\label{prop non-trivial}
Let $k$ be a CM-field, $p$ an odd prime number that splits completely in $k$, and $\Sigma$ a $p$-adic CM-type of $k$. If $[k:\Q]\geq 4$, then $X_\Sigma(\widetilde{k})$ is not trivial. 
\end{prop}
Since $X(\wtk)$ is a quotient of $X_\Sigma(\wtk)$, the non-triviality of $X_\Sigma(\wtk)$ follows from Proposition~\ref{prop Fujii non-trivial} in the case where $[k:\Q]\geq 8$. 

We show Proposition~\ref{prop non-trivial} in two ways. The first way is to use a Tate sequence and the second way is to use Hochschild--Serre spectral sequence. 

In the present article, we provide some consequences of Proposition~\ref{prop non-trivial}.

\subsection{Non-coincidence of unramified and $\Sigma$-ramified Iwasawa modules}
In the present article, we are interested in the following weak form of GGC. 
\begin{conjecture}[a weak form of Greenberg's Generalized Conjecture (wGGC) {\cite[GC(2)]{WingbergPre}}]\label{conjecture wGGC}
Let $k/\Q$ be a finite extension and $p$ a prime number. If $X(\wtk)$ is not trivial, then $X(\wtk)$ has a non-trivial pseudo-null $\Z_p\bbr{\Gal(\wtk/k)}$-submodule. 
\end{conjecture}
Concerning wGGC, Fujii showed the following criterion. 

\begin{thm}[{\cite[Theorem~2]{Fujii2014}}]\label{thm Fujii criterion}
Let $k$ be an imaginary quadratic field, $p$ an odd prime number that splits in $k$ as $p=\pe\overline{\pe}$. Then the following are equivalent. 
\begin{enumerate}[label=(\roman*)]
\item $X(\widetilde{k})$ has a non-trivial pseudo-null $\Z_p\bbr{\Gal(\widetilde{k}/k)}$-submodule. \label{Fujii 1}

\item $L(\widetilde{k})\neq M_\pe(\widetilde{k})$ holds. \label{Fujii 2}
\end{enumerate}
\end{thm}
Theorem~\ref{thm Fujii criterion} can be seen as an analogue of Ozaki's result in \cite{Ozaki1997}, which deals with the cyclotomic $\Z_p$-extension over a totally real number field (see \cite[Theorem~A]{Fujii2014} as well). We mention that Murakami \cite{Murakami2023} used Theorem~\ref{thm Fujii criterion} to show a sufficient condition for wGGC for imaginary quadratic fields. 

It is an interesting problem to extend Theorem~\ref{thm Fujii criterion} to higher degree CM-fields. From this point of view, formulating the following conjecture may be of some significance. 
\begin{conjecture}\label{conjecture different}
Let $k$ be a CM-field, $p$ an odd prime number that splits completely in $k$, and $\Sigma$ a $p$-adic CM-type of $k$. If $[k:\Q]\geq 4$, then $L(\wtk)\neq M_\Sigma(\wtk)$ holds. 
\end{conjecture}
It can be easily shown by Proposition~\ref{prop non-trivial} that wGGC for $(k, p)$ implies Conjecture~\ref{conjecture different} (see Proposition~\ref{prop wGGC different}).

\subsection{Non-freeness Conjecture and Pseudo-null Iwasawa modules}
We investigate a relation between Conjecture~\ref{conjecture different} and a conjecture called Non-freeness Conjecture. 

For any algebraic extension $K/\Q$ and a prime number $p$, we let $\calL(K)/K$ denote the maximal unramified pro-$p$-extension and $\calG(K)\coloneqq\Gal(\calL(K)/K)$ its Galois group. We also let $\calL^{(2)}(K)/K$ denote the maximal unramified metabelian pro-$p$-extension and put $\calG^{(2)}(K)\coloneqq\Gal(\calL^{(2)}(K)/K)$. Then $X(K)$ (resp. $\calG^{(2)}(K)$) is the abelian (resp. metabelian) quotient of $\calG(K)$. 

 For a finite extension $k/\Q$ and a prime number $p$, we let $k^\cy_\infty/k$ denote the cyclotomic $\Z_p$-extension. The following conjecture is studied in non-abelian Iwasawa theory in the sense of Ozaki \cite{Ozaki2007}.
\begin{conjecture}[Non-freeness Conjecture]\label{conjecture NF}
Let $k/\Q$ be a finite extension and $p$ a prime number. Then $\calG(k^\cy_\infty)$ is not a non-abelian free pro-$p$-group. 
\end{conjecture}
This conjecture was formulated by Ozaki (see \cite[p.531, Conjecture]{Fujii2022}). Non-freeness Conjecture is related to wGGC as the following result shows. 
\begin{thm}[{\cite[Proposition~2.1]{Fujii2011}}]\label{thm wGGC NF}
Let $k/\Q$ be a finite extension and $p$ a prime number that  splits completely in $k$. Assume that wGGC is true for $(k, p)$. Then $\calG^{(2)}(k^\cy_\infty)$ is not a non-abelian free metabelian pro-$p$-group and, in particular, Non-freeness Conjecture is true for $(k, p)$. 
\end{thm}
While \cite[Proposition~2.1]{Fujii2011} is stated in terms of derived depth, the statement in terms of free metabelian pro-$p$-group can be found in the second remark after \cite[Theorem~1.1]{Fujii2011RIMS}. 

We show a variant of Theorem~\ref{thm wGGC NF} in the case of CM-fields. To state the result, we introduce a weak form of original Greenberg's Conjecture. 
\begin{conjecture}[A weak form of Greenberg's Conjecture (wGC)]\label{conjecture wGC}
Let $k/\Q$ be a totally real finite extension and $p$ a prime number. If $X(k^\cy_\infty)$ is not trivial, then $X(k^\cy_\infty)$ has a non-trivial finite $\Z_p\bbr{\Gal(k^\cy_\infty/k)}$-submodule. 
\end{conjecture}
Conjecture~\ref{conjecture wGC} was proposed by Nguyen Quang Do \cite[(GF)]{Nguyen2006} (the above formulation can be found in \cite[(WGC)]{Nguyen2017}). See \cite{Greenberg1976} for original Greenberg's Conjecture. 
Note that, for a totally real finite extension $k/\Q$, Conjecture~\ref{conjecture wGC} coincides with wGGC, under $p$-adic Leopoldt's Conjecture. Indeed, the extensions $\wtk/k$ and $k^\cy_\infty/k$ coincide.

\subsection{Main results}
In the present article, we show the following result, which shows a relation between Conjecture~\ref{conjecture different}, wGC, and Non-freeness Conjecture. 
\begin{thm}\label{thm NC wGC to NF}
Let $k$ be a CM-field of degree greater than $2$, $p$ an odd prime number that splits completely in $k$, and $\Sigma$ a $p$-adic CM-type of $k$. Assume that Conjecture~\ref{conjecture different} is true for $(k, \Sigma)$, and wGC is true for $(k^+, p)$. Then $\calG^{(2)}(k^\cy_\infty)$ is not a free metabelian pro-$p$-group and, in particular, Non-freeness Conjecture is true for $(k, p)$. 
\end{thm}
The strategy for Theorem~\ref{thm NC wGC to NF} is as follows. To show the contraposition, suppose that wGC is true for $(k^+, p)$, and $\calG^{(2)}(k^\cy_\infty)$ is a free metabelian pro-$p$-group. We first observe that $L(k^\cy_\infty)=M_\Sigma(k^\cy_\infty)$ holds. Then we show that the identity $X_\Sigma(k^\cy_\infty)\to X(k^\cy_\infty)$ lifts to the identity $\Gal(M_\Sigma(\wtk)/k^\cy_\infty)\to \Gal(L(\wtk)/k^\cy_\infty)$ by using the freeness of $\calG^{(2)}(k^\cy_\infty)$, so we obtain $L(\wtk)=M_\Sigma(\wtk)$. 

If we focus on the case of quartic CM-fields with a certain assumption, we obtain the following result, which gives an analogue of Theorem~\ref{thm Fujii criterion} and shows the inverse implication of Theorem~\ref{thm NC wGC to NF}. 
\begin{thm}\label{thm main}
Let $k$ be a quartic CM-field, $p$ an odd prime number that splits completely in $k$, and $\Sigma$ a $p$-adic CM-type of $k$. Suppose that $X(k^\cy_\infty)\simeq \Z_p^2$ as a $\Z_p\bbr{\Gal(k^\cy_\infty/k)}$-module. Then the following are equivalent. 
\begin{enumerate}[label=(\roman*)]
\item $X(\widetilde{k})$ is pseudo-null as a $\Z_p\bbr{\Gal(\widetilde{k}/k)}$-module (GGC). \label{main item 1}

\item If $X(\wtk)$ is not trivial, then $X(\widetilde{k})$ has a non-trivial pseudo-null $\Z_p\bbr{\Gal(\widetilde{k}/k)}$-submodule (wGGC).\label{main item 2}

\item $L(\widetilde{k})\neq M_\Sigma(\widetilde{k})$ holds (Conjecture~\ref{conjecture different}). \label{main item 3}
\item $\calG^{(2)}(k^\cy_\infty)$ is not a free metabelian pro-$p$-group (in particular, Non-freeness Conjecture follows). \label{main item 4}
\end{enumerate}
\end{thm}
In \S\ref{section criterion}, we provide some examples of GGC by using Theorem~\ref{thm main}. We show a sufficient condition for Condition~\ref{main item 4} in terms of the ideal class groups of certain finite extensions of $k$ (Proposition~\ref{prop criterion}). 
\begin{rem}\label{rem Bandini}
In the setting of Theorem~\ref{thm main}, if we suppose that $k$ is an imaginary biquadratic field, then the validity of GGC is shown by Bandini \cite[Theorem~9 and Corollary~13]{Bandini2007}. 
\end{rem}

\subsection*{Structure of the present article}
The structure of the present article is as follows. In \S\ref{section preliminaries}, we gather some basic notation and algebraic lemmas. In \S\ref{section Lemmas on CM-fields}, we show some basic properties of Iwasawa modules over CM-fields. We prove Proposition~\ref{prop non-trivial} in \S\ref{section Non-triviality}. In \S\ref{section general CM-fields}, we show some results, including Theorem~\ref{thm NC wGC to NF}, that do not place assumptions on the degree of CM-fields. In \S\ref{section proof main}, we prove Theorem~\ref{thm main}. In \S\ref{section criterion}, we show a criterion for GGC and provide concrete examples.

\section{Preliminaries}\label{section preliminaries}
\subsection{Modules over completed group algebras}\label{subsection pseudo-null}
Let $p$ be a prime number, $\Gamma$ a pro-$p$-group isomorphic to $\Z_p^d$ for some $d\geq 0$, and put $\Lambda\coloneqq\Z_p\bbr{\Gamma}$. For a $\Lambda$-module $M$, we let $M_\Gamma$ denote the $\Gamma$-coinvariant of $M$. 

For any profinite ring $A$ and a left profinite $A$-module $M$, we let $\pd_AM$ denote the projective dimension of $M$ in the abelian category of profinite $A$-modules. Note that, for a finitely generated $\Lambda$-module $M$, the topological projective dimension $\pd_\Lambda M$ coincides with the projective dimension of $M$ as an abstract $\Lambda$-module. This is because $\Lambda$ is noetherian. For any profinite group $G$, we let $I_p(G)$ denote the $p$-adic augmentation kernel of $G$, that is, the kernel of the augmentation map $\Z_p\bbr{G}\to \Z_p$. We also let $\cd_pG$ denote the $p$-cohomological dimension of $G$. It is known that $\pd_{\Z_p\bbr{G}}\Z_p=\cd_pG$ holds (see \cite[(5.2.12)~Corollary]{NSW2008}). 

Let $M$ be a finitely generated $\Lambda$-module. We say that $M$ is pseudo-null as a $\Lambda$-module if the height of the annihilator ideal $\Ann_\Lambda M$ is greater than $1$. It is known that $M$ is pseudo-null if and only if there exist two relatively prime annihilators $f, g\in\Ann_\Lambda M$. Here we note that $\Lambda$ is a unique factorization domain. 

\begin{lem}\label{lem tor pn}
Let $H$ be a closed subgroup of $\Gamma$ such that $\Gamma/H\simeq\Z_p$ and $M$ a $\Lambda$-module that is finitely generated as a $\Z_p\bbr{H}$-module. Then the following are equivalent. 
\begin{enumerate}[label=(\roman*)]
\item $M$ is pseudo-null as a $\Lambda$-module. \label{item pn}

\item $M$ is torsion as a $\Z_p\bbr{H}$-module. \label{item tors}
\end{enumerate}
\end{lem}
\begin{proof}
This lemma is classically known. A proof can be found, for example, in \cite[Lemma~2.3]{Fujii2011}.  Here we briefly show the implication \ref{item tors}$\Rightarrow$\ref{item pn}, which we only use in the present article. 

Take a lift $\gamma\in\Gamma$ of a topological generator of $\Gamma/H$. Then there is an isomorphism $\Lambda\simeq\Z_p\bbr{H}\bbr{T}$ defined by $\gamma\leftrightarrow 1+T$. Consider a $\Z_p\bbr{H}$-linear endomorphism of $M$ defined by multiplication by $T$. By Cayley--Hamilton theorem, there is a monic polynomial $f(T)\in\Z_p\bbr{H}\br{T}$ that annihilates $M$. Now suppose that $M$ is torsion as a $\Z_p\bbr{H}$-module. Then there is a non-zero annihilator $g\in\Ann_{\Z_p\bbr{H}}M$. Since $f(T)$ and $g$ are relatively prime, we conclude that $M$ is pseudo-null as a $\Lambda$-module. 
\end{proof}

\subsection{A lemma on free pro-$p$-groups}
Let $p$ be a prime number and $\calG$ a pro-$p$-group. We let $D_1(\calG)\coloneqq\overline{[\calG, \calG]}$ denote the closed commutator subgroup of $\calG$ and $\calG_\ab\coloneqq \calG/D_1(\calG)$ the abelianization of $\calG$. Put $D_2(\calG)\coloneqq D_1(D_1(\calG))$. 
\begin{lem}\label{lem rank 2 free}
Let $F$ be a free pro-$p$-group of rank $2$. Then $D_1(F)_\ab$ is a free $\Z_p\bbr{F_\ab}$-module of rank $1$. 
\end{lem}
\begin{proof}
Lyndon's resolution attached to $F_\ab$ gives an exact sequence
\[
\begin{tikzpicture}
\node (a) at (0,0) {$0$}; 
\node (b) at (2.2,0) {$D_1(F)_\ab$}; 
\node (c) at (5,0) {$\Z_p\bbr{F_\ab}^2$}; 
\node (d) at  (7.8,0) {$\Z_p\bbr{F_\ab}$}; 
\node (e) at (10.2,0) {$\Z_p$}; 
\node (f) at (12,0) {$0$};  

\draw[->] (a) -- (b); \draw[->] (b) -- (c); \draw[->] (c) -- (d); \draw[->] (d) -- (e); \draw[->] (e) -- (f); 
\end{tikzpicture}
\]
of $\Z_p\bbr{F_\ab}$-modules (see \cite[Proposition~1.1]{Nguyen1984}). So the $\Z_p\bbr{F_\ab}$-rank of $D_1(F)_\ab$ is $1$. Since $F_\ab$ is isomorphic to $\Z_p^2$, the projective dimension of $\Z_p$ as a $\Z_p\bbr{F_\ab}$-module is $2$. Hence $D_1(F)_\ab$ is projective as a $\Z_p\bbr{F_\ab}$-module by the above exact sequence. Since $\Z_p\bbr{F_\ab}$ is a local ring, it follows that $D_1(F)_\ab$ is a free $\Z_p\bbr{F_\ab}$-module. 
\end{proof}

\section{Lemmas on CM-fields}\label{section Lemmas on CM-fields}
In this section, we state some lemmas concerning CM-fields. The facts stated in this section are not new, but we provide proofs of some lemmas for completeness. 

Let $k$ be a CM-field of  degree $2d$ with $d\geq 1$, $p$ an odd prime number that splits completely in $k$, and $\Sigma$ a $p$-adic CM-type of $k$. Recall that a set $\Sigma$ of $p$-adic primes of $k$ is called a $p$-adic CM-type if it satisfies $S_p(k)=\Sigma\cup\overline{\Sigma}$ and $\Sigma\cap\overline{\Sigma}=\emptyset$. Here we write $\overline{\Sigma}=\{\overline{\pe}\mid \pe\in\Sigma\}$, where $\overline{\pe}$ denotes the complex conjugate of $\pe$. Put $\Gamma\coloneqq\Gal(\wtk/k)$ and $\Lambda\coloneqq\Z_p\bbr{\Gamma}$. We let $k^+_\infty/k^+$ denote the cyclotomic $\Z_p$-extension. 
\begin{lem}[{\cite[Lemma~2.10]{Yanagisawa2026}}]\label{lem LM}
Suppose that $X(k^+_\infty)$ is trivial. Then $L(k^\cy_\infty)=M_\Sigma(k^\cy_\infty)$ holds. 
\end{lem}
\begin{proof}
We have $X_p(k^+_\infty)=0$ by \cite[Corollary~1]{Ozaki1997}. Thus the generator of $\Gal(k^\cy_\infty/k^+_\infty)$ acts as inverse on $X_p(k^\cy_\infty)$. Hence every intermediate field of $M_p(k^\cy_\infty)/k^\cy_\infty$ is a Galois extension over $k^+_\infty$ and, in particular, the extension $M_\Sigma(k^\cy_\infty)/k^+_\infty$ is a Galois extension. Then, by considering the action of the generator of $\Gal(k^\cy_\infty/k^+_\infty)$, we see that $M_\Sigma(k^\cy_\infty)/k^\cy_\infty$ is $\overline{\Sigma}$-ramified. Therefore, we conclude that $M_\Sigma(k^\cy_\infty)/k^\cy_\infty$ is unramified. This means that $M_\Sigma(k^\cy_\infty)\subseteq L(k^\cy_\infty)$ holds. The reverse inclusion is clear. 
\end{proof}
\begin{lem}\label{lem four-term}
Suppose that $X(k^+_\infty)$ is trivial. Then there is an exact sequence
\[
\begin{tikzpicture}
\node (a) at (0,0) {$\bigwedge^2\Gal(\wtk/k^\cy_\infty)$}; 
\node (b) at (3.7,0) {$X_\Sigma(\wtk)_{\Gal(\wtk/k^\cy_\infty)}$}; 
\node (c) at (6.9,0) {$X(k^\cy_\infty)$}; 
\node (d) at  (9.6,0) {$\Gal(\wtk/k^\cy_\infty)$}; 
\node (e) at (11.5,0) {$0$}; 

\draw[->] (a) -- (b); \draw[->] (b) -- (c); \draw[->] (c) -- (d); \draw[->] (d) -- (e); 
\end{tikzpicture}
\]
of $\Z_p\bbr{\Gal(k^\cy_\infty/k)}$-modules, where $\bigwedge^2 (-)$ denotes the exterior product over $\Z_p$. 
\end{lem}
\begin{proof}
By \cite[Proposition~1]{Ozaki1997}, the extension $\wtk/k^\cy_\infty$ is unramified. Thus the maximal subfield of $M_\Sigma(\wtk)$ that is abelian over $k^\cy_\infty$ coincides with $M_\Sigma(k^\cy_\infty)$, which coincides with $L(k^\cy_\infty)$ by Lemma~\ref{lem LM}. Hence, by the five-term exact sequence in group homology, we obtain the desired exact sequence (for details, see \cite[Lemma~2.5]{Yanagisawa2026}, for example). 
\end{proof}

\begin{lem}\label{lem fg tors}
$X_\Sigma(\wtk)$ is finitely generated and torsion as a $\Lambda$-module. 
\end{lem}
\begin{proof}
See \cite[Theorem~1.2.2]{HT1994}.  
\end{proof}
The following lemma was established by Fujii in the proof of \cite[Proposition~1]{Fujii2014} in the case where $k$ is an imaginary quadratic field, and here we generalize it to CM-fields. 
\begin{lem}\label{lem exact}
There is an exact sequence
\[
\begin{tikzpicture}
\node (a) at (0,0) {$0$}; 
\node (b) at (1.7,0) {$\Lambda^d$}; 
\node (c) at (3.7,0) {$X_p(\wtk)$}; 
\node (d) at  (6,0) {$X_\Sigma(\wtk)$}; 
\node (e) at (7.8,0) {$0$}; 

\draw[->] (a) -- (b); \draw[->] (b) -- (c); \draw[->] (c) -- (d); \draw[->] (d) -- (e); 
\end{tikzpicture}
\]
of $\Lambda$-modules, where the map $X_p(\wtk)\to X_\Sigma(\wtk)$ is the restriction. 
\end{lem}
\begin{proof}
Put $Y\coloneqq\Gal(M_p(\wtk)/M_\Sigma(\wtk))=\Ker(X_p(\wtk)\to X_\Sigma(\wtk))$. For each $\Pe\in S_p(\wtk)$, we let $I_\Pe$ denote the inertia group of $\Pe$ in $X_p(\wtk)$. Then there is a natural surjection
\[
\bigoplus_{\pe\in\overline{\Sigma}}\tgen{I_\Pe\mid \Pe\in S_\pe(\wtk)}\to Y
\]
of $\Lambda$-modules. Here we note that each summand of the direct sum is a $\Lambda$-module. 

Take $\pe\in\overline{\Sigma}$, $\Pe_0\in S_\pe(\wtk)$, and let $Z_\pe$ be the decomposition field of $\pe$ in $\wtk/k$. Put $\pe_0\coloneqq\Pe_0\cap Z_\pe$. By local class field theory, the $\Z_p$-rank of the inertia group of $\pe$ in $\wtk/k$ is at most $1$. On the other hand, since $\wtk$ contains $k^\cy_\infty$, we see that the $\Z_p$-rank is equal to $1$. Hence it follows that the $\Z_p$-rank of the decomposition group $\Gal(\wtk/Z_\pe)$ is $1$ or $2$. Now we claim that there is a surjection $\Z_p\bbr{\Gal(\wtk/Z_\pe)}\to I_{\Pe_0}$. Indeed, if $\rank_{\Z_p}\Gal(\wtk/Z_\pe)=1$, then it follows from Iwasawa's theorem \cite[Theorem~25~(i)]{Iwasawa1973}, and if $\rank_{\Z_p}\Gal(\wtk/Z_\pe)=2$, then it follows from Wintenberger's theorem \cite[Theorem~(i)]{Wintenberger1980}. Thus, by operating the completed tensor product $\Lambda\mathrel{\widehat{\otimes}}_{\Z_p\bbr{\Gal(\wtk/Z_\pe)}}(-)$ to the surjection, we obtain a surjection
\[
\Lambda\to\tgen{I_\Pe\mid \Pe\in S_\pe(\wtk)}
\]
of $\Lambda$-modules. By gathering the surjections for all $\pe\in\overline{\Sigma}$, we see that there is a surjection $\Lambda^d\to Y$. So there is an exact sequence
\[
\begin{tikzpicture}
\node (b) at (1.7,0) {$\Lambda^d$}; 
\node (c) at (3.7,0) {$X_p(\wtk)$}; 
\node (d) at  (6,0) {$X_\Sigma(\wtk)$}; 
\node (e) at (7.8,0) {$0$}; 

\draw[->] (b) -- (c); \draw[->] (c) -- (d); \draw[->] (d) -- (e); 
\end{tikzpicture}
\]
of $\Lambda$-modules. Now $X_\Sigma(\wtk)$ is torsion as a $\Lambda$-module by Lemma~\ref{lem fg tors}, and the $\Lambda$-rank of $X_p(\wtk)$ is $d$ by Greenberg's result \cite[Proposition~2]{Greenberg1978}. Therefore, the kernel of the map $\Lambda^d\to X_p(\wtk)$ is torsion as a $\Lambda$-module, and hence it is trivial. This completes the proof. 
\end{proof}
\begin{rem}
As mentioned in \cite[Theorem~D, Remark]{Fujii2011}, the fact where we referred to Wintenberger's theorem \cite[Theorem~(i)]{Wintenberger1980} can alternatively be shown by the fact that the Galois group of the maximal pro-$p$-extension of $\Q_p$ is a free pro-$p$-group of rank $2$ (see \cite[(7.5.11)~Theorem]{NSW2008}). Indeed, if we let $\widetilde{\Q_p}/\Q_p$ denote the maximal abelian pro-$p$-extension, then it is the $\Z_p^2$-extension over $\Q_p$, and the Galois group of the maximal abelian pro-$p$-extension over $\widetilde{\Q_p}$ is a free $\Z_p\bbr{\Gal(\widetilde{\Q_p}/\Q_p)}$-module of rank $1$ by Lemma~\ref{lem rank 2 free}. Since there is no non-trivial unramified pro-$p$-extension over $\widetilde{\Q_p}$, the assertion follows. 
\end{rem}

\begin{lem}\label{lem pn-free}
The module $X_\Sigma(\widetilde{k})$ has no non-trivial pseudo-null $\Lambda$-submodule. 
\end{lem}
\begin{proof}
One can refer to \cite[Proposition~3.5]{HO2010} for the algebraic arguments below. But we give a proof following Fujii's proof of \cite[Proposition~1]{Fujii2014}, because our setting is somewhat more elementary than the non-commutative setting in \cite{HO2010}. 

Let $B$ be a pseudo-null $\Lambda$-submodule of $X_\Sigma(\wtk)$. It suffices to show the triviality of $B$. Let $C$ be the preimage of $B$ with respect to the surjection $X_p(\wtk)\to X_\Sigma(\wtk)$. Then, by Lemma~\ref{lem exact}, we obtain an exact sequence
\[
\begin{tikzpicture}
\node (a) at (0,0) {$0$}; 
\node (b) at (1.7,0) {$\Lambda^d$}; 
\node (c) at (3.7,0) {$C$}; 
\node (d) at  (5.7,0) {$B$}; 
\node (e) at (7.6,0) {$0$}; 

\draw[->] (a) -- (b); \draw[->] (b) -- (c); \draw[->] (c) -- (d); \draw[->] (d) -- (e); 
\end{tikzpicture}
\]
of $\Lambda$-modules. Observe that $\Tor_\Lambda C$ injects into $B$. Since $B$ is pseudo-null, it follows that $\Tor_\Lambda C$ is a pseudo-null submodule of $X_p(\wtk)$. On the other hand, it is shown in \cite[Corollary~3.2]{Nguyen1984} that $X_p(\wtk)$ has no non-trivial pseudo-null $\Lambda$-submodule, so $\Tor_\Lambda C$ is trivial. 

Since $B$ is pseudo-null, we can take two relatively prime annihilators $f, g\in \Ann_\Lambda B$. Consider an exact commutative diagram
\[
\begin{tikzpicture}
\node (a) at (0,0) {$0$}; 
\node (b) at (1.7,0) {$\Lambda^d$}; 
\node (c) at (3.7,0) {$C$}; 
\node (d) at  (5.7,0) {$B$}; 
\node (e) at (7.6,0) {$0$}; 

\node (a1) at (0,2) {$0$}; 
\node (b1) at (1.7,2) {$\Lambda^d$}; 
\node (c1) at (3.7,2) {$C$}; 
\node (d1) at  (5.7,2) {$B$}; 
\node (e1) at (7.6,2) {$0$};

\draw[->] (a) -- (b); \draw[->] (b) -- (c); \draw[->] (c) -- (d); \draw[->] (d) -- (e); 

\draw[->] (a1) -- (b1); \draw[->] (b1) -- (c1); \draw[->] (c1) -- (d1); \draw[->] (d1) -- (e1); 

\draw[->] (b1) -- node[auto][swap]{$\scriptstyle \times g$} (b); \draw[->] (c1) -- node[auto][swap]{$\scriptstyle \times g$} (c); \draw[->] (d1) -- node[auto][swap]{$\scriptstyle \times g$} (d); 
\end{tikzpicture}
\]
of $\Lambda$-modules. Then, by Snake lemma, we obtain an exact sequence
\[
\begin{tikzpicture}
\node (a) at (0,0) {$C\br{g}$}; 
\node (b) at (2,0) {$B\br{g}$}; 
\node (c) at (4.3,0) {$(\Lambda/(g))^d$}; 

\draw[->] (a) -- (b); \draw[->] (b) -- (c); 
\end{tikzpicture}
\]
of $\Lambda$-modules. Here, we note that $C\br{g}$ is trivial since $\Tor_\Lambda C$ is trivial, and $B\br{g}=B$ because $g$ annihilates $B$. Hence we obtain an injection $B\to (\Lambda/(g))^d$. By taking the $f$-torsion parts, we obtain an injection $B\br{f}\to (\Lambda/(g))^d\br{f}$. Here, we have $B\br{f}=B$ because $f$ annihilates $B$, and $(\Lambda/(g))^d\br{f}$ is trivial because $f$ and $g$ are relatively prime. Therefore, it follows that $B$ is trivial. This completes the proof. 
\end{proof}

\section{Non-triviality of Iwasawa modules over CM-fields}\label{section Non-triviality}
We keep the setting in \S\ref{section Lemmas on CM-fields}. 
\subsection{First proof of Proposition~\ref{prop non-trivial}}\label{subsection first proof}
\begin{lem}\label{lem pd d-1}
An inequality $\pd_\Lambda X_p(\wtk)\geq d-1$ holds. 
\end{lem}
This lemma follows from \cite[Proposition~3.4]{Nguyen1984}, but here we recall the proof for its importance. 

Let $G$ (resp. $H$) be the Galois group of the maximal $p$-ramified extension over $k$ (resp. $\wtk$), and put $R\coloneqq\Z_p\bbr{G}$. 

The following exact sequence is known as Tate sequence. It is the right-hand vertical sequence of \cite[(5.6.7)~Proposition]{NSW2008}. 
\begin{lem}\label{lem Tate sequence}
There is an exact sequence
\[
\begin{tikzpicture}
\node (a) at (0,0) {$0$}; 
\node (b) at (1.7,0) {$X_p(\wtk)$}; 
\node (c) at (4.2,0) {$\Lambda\otimes_R I_p(G)$}; 
\node (d) at  (6.3,0) {$\Lambda$}; 
\node (e) at (7.7,0) {$\Z_p$}; 
\node (f) at (9,0) {$0$};

\draw[->] (a) -- (b); \draw[->] (b) -- (c); \draw[->] (c) -- (d); \draw[->] (d) -- (e); \draw[->] (e) -- (f); 
\end{tikzpicture}
\]
of $\Lambda$-modules. 
\end{lem}
\begin{proof}
Take the $H$-homology long exact sequence of $0\to I_p(G)\to R\to\Z_p\to 0$. 
\end{proof}
\begin{lem}\label{lem pd 1}
An inequality $\pd_\Lambda(\Lambda\otimes_R I_p(G))\leq 1$ holds.
\end{lem}
\begin{proof}
By \cite[Proposition~(8.3.18)]{NSW2008}, we have $\cd_pG\leq 2$, so $\pd_R\Z_p\leq 2$. Hence we have $\pd_R I_p(G)\leq 1$ by a dimension shifting argument. Take a projective resolution
\[
\begin{tikzpicture}
\node (a) at (0,0) {$0$}; 
\node (b) at (1.6,0) {$P_1$}; 
\node (c) at (3.6,0) {$P_0$}; 
\node (d) at  (5.8,0) {$I_p(G)$}; 
\node (e) at (7.6,0) {$0$}; 

\draw[->] (a) -- (b); \draw[->] (b) -- (c); \draw[->] (c) -- (d); \draw[->] (d) -- (e); 
\end{tikzpicture}
\]
of $I_p(G)$. Then, by taking the $H$-homology long exact sequence, we obtain an exact sequence
\[
\begin{tikzpicture}
\node (a) at (0.3,0) {$\Homo_1(H, I_p(G))$}; 
\node (b) at (3,0) {$\Lambda\otimes_R P_1$}; 
\node (c) at (5.7,0) {$\Lambda\otimes_RP_0$}; 
\node (d) at  (8.3,0) {$\Lambda\otimes_RI_p(G)$}; 
\node (e) at (10,0) {$0$}; 

\draw[->] (a) -- (b); \draw[->] (b) -- (c); \draw[->] (c) -- (d); \draw[->] (d) -- (e); 
\end{tikzpicture}
\]
of $\Lambda$-modules. To show the inequality $\pd_\Lambda(\Lambda\otimes_R I_p(G))\leq 1$, it suffices to show the triviality of $\Homo_1(H, I_p(G))$. Since we have $\Homo_i(H, R)=0$ for $i=1, 2$, there is an isomorphism $\Homo_1(H, I_p(G))\simeq\Homo_2(H, \Z_p)$. Now $\Homo_2(H, \Z_p)$ is trivial by weak Leopoldt's Conjecture (see \cite[Theorem~2.2]{Nguyen1984}), so we conclude that $\Homo_1(H, I_p(G))$ is trivial, as desired. 
\end{proof}
The proof of Lemma~\ref{lem pd d-1} can be completed by combining the above lemmas and an inequality $\pd_\Lambda\Z_p=\cd_p\Gamma\geq d+1$. 
Now it is ready to prove Proposition~\ref{prop non-trivial}. 
\begin{proof}[Proof of Proposition~\ref{prop non-trivial}]
To show the contraposition, we assume that $X_\Sigma(\wtk)$ is trivial. Then we can use Lemma~\ref{lem exact} to obtain an isomorphism $X_p(\wtk)\simeq \Lambda^d$. By Lemma~\ref{lem pd d-1}, we obtain an inequality $0\geq d-1$, so $d=1$. This completes the proof. 
\end{proof}

\subsection{Second proof of Proposition~\ref{prop non-trivial}}
Here we give another proof of Proposition~\ref{prop non-trivial}. Keep the notation in \S\ref{subsection first proof}. 

Since $\Homo_2(H, \Z_p)$ is trivial by weak Leopoldt's Conjecture (see \cite[Theorem~2.2]{Nguyen1984}), we obtain an exact sequence
\[
\begin{tikzpicture}
\node (a) at (0,0) {$\Homo_3(G, \Z_p)$}; 
\node (b) at (2.5,0) {$\Homo_3(\Gamma, \Z_p)$}; 
\node (c) at (5.5,0) {$\Homo_1(\Gamma, \Homo_1(H,\Z_p))$}; 

\draw[->] (a) -- (b); \draw[->] (b) -- (c); 
\end{tikzpicture}
\]
of $\Z_p$-modules by Hochschild--Serre spectral sequence (see \cite[Theorem~7.2.4]{RZ2010}, for example). Here, $\Homo_3(G,\Z_p)$ is trivial, because $\cd_p G\leq 2$ by \cite[Proposition~(8.3.18)]{NSW2008}. Thus the above exact sequence gives an injection
\[
\bigwedge^3\Gamma\to\Homo_1(\Gamma, X_p(\wtk))
\]
of $\Z_p$-modules, where $\bigwedge^3\Gamma$ denotes the exterior power over $\Z_p$. Now, if we suppose that $X_\Sigma(\wtk)$ is trivial, then we obtain $X_p(\wtk)\simeq \Lambda^d$ by Lemma~\ref{lem exact}, so $\Homo_1(\Gamma, X_p(\wtk))$ is trivial. Since the $\Z_p$-rank of $\bigwedge^3\Gamma$ is at least $(d+1)d(d-1)/6$, we obtain $d=1$. This completes the proof.

\subsection{Non-triviality of the unramified Iwasawa module}\label{subsection unramified}
As we remarked in Remark~\ref{rem Leopoldt}, it is not necessary to assume the validity of $p$-adic Leopoldt's Conjecture in Proposition~\ref{prop Fujii non-trivial}. Here we briefly recall the proof of Proposition~\ref{prop Fujii non-trivial} for convenience. 
Let $k/\Q$ be a totally imaginary finite extension, $d$ the number of complex places of $k$, and $p$ a prime number that  splits completely in $k$. Let $k^\cy_n$ denote the $n$-th layer of $k^\cy_\infty/k$ and put $E\coloneqq\varprojlim_n (\Z_p\otimes_{\Z}E(k^\cy_n))$, where $E(k^\cy_n)$ denotes the unit group of $k^\cy_n$. By \cite[Lemma~3.7]{Fujii2022}, the $\Z_p$-rank of $E_{\Gal(k^\cy_\infty/k)}$ is $d$. 

Now suppose that $X(\wtk)$ is trivial. Then $\calG(\wtk)$ is trivial as well. Since $\wtk/k^\cy_\infty$ is unramified by \cite[Lemma~3.1]{Fujii2022}, it follows that $\calG(k^\cy_\infty)=\Gal(\wtk/k^\cy_\infty)$. Then, by \cite[Theorem~C]{Ozaki2007} for $i=2$, there is a surjection
\[
E\to \Homo_2(\Gal(\wtk/k^\cy_\infty), \Z_p)\simeq \bigwedge^2\Gal(\wtk/k^\cy_\infty)
\]
of $\Z_p\bbr{\Gal(k^\cy_\infty/k)}$-modules. Since $\Gal(k^\cy_\infty/k)$ acts trivially on the right-hand side, the surjection factors through $E_{\Gal(k^\cy_\infty/k)}$. The $\Z_p$-rank of $\bigwedge^2\Gal(\wtk/k^\cy_\infty)$ is at least $d(d-1)/2$, so we obtain an inequality
\[
d\geq\frac{d(d-1)}{2}, 
\]
by comparing the $\Z_p$-ranks. Hence we obtain $d\leq 3$. This completes the proof.

\section{Results on general CM-fields}\label{section general CM-fields}
\subsection{Non-coincidence of unramified and $\Sigma$-ramified Iwasawa modules}
\begin{prop}\label{prop wGGC different}
Let $k$ be a CM-field of degree greater than $2$, $p$ an odd prime number that splits completely in $k$, and $\Sigma$ a $p$-adic CM-type of $k$. If wGGC is true for $(k, p)$, then Conjecture~\ref{conjecture different} is true for $(k, \Sigma)$. 
\end{prop}
\begin{proof}
 If $X(\wtk)$ is trivial, then the assertion follows from Proposition~\ref{prop non-trivial}. On the other hand, if $X(\wtk)$ is not trivial, then the assertion follows from Lemma~\ref{lem pn-free}. 
\end{proof}

\begin{thm}\label{thm general}
Let $k$ be a CM-field of degree greater than $2$, $p$ an odd prime number that splits completely in $k$, and $\Sigma$ a $p$-adic CM-type of $k$. Suppose that the characteristic ideal $\chara_\Lambda X_\Sigma(\wtk)$ is a prime ideal. Then the following are equivalent. 
\begin{enumerate}[label=(\roman*)]
\item $X(\widetilde{k})$ is pseudo-null as a $\Z_p\bbr{\Gal(\widetilde{k}/k)}$-module (GGC). 

\item If $X(\wtk)$ is not trivial, then $X(\widetilde{k})$ has a non-trivial pseudo-null $\Z_p\bbr{\Gal(\widetilde{k}/k)}$-submodule (wGGC). 

\item $L(\widetilde{k})\neq M_\Sigma(\widetilde{k})$ holds (Conjecture~\ref{conjecture different}). 
\end{enumerate}
\end{thm}
\begin{proof}
The implication \ref{main item 1}$\Rightarrow$\ref{main item 2} is clear. The implication \ref{main item 2}$\Rightarrow$\ref{main item 3} is Proposition~\ref{prop wGGC different}. 

From now on, we show the implication \ref{main item 3}$\Rightarrow$\ref{main item 1}. Put $N\coloneqq\Ker(X_\Sigma(\wtk)\to X(\wtk))=\Gal(M_\Sigma(\wtk)/L(\wtk))$. Then $N$ is a non-trivial $\Lambda$-submodule of $X_\Sigma(\wtk)$, so $N$ is not pseudo-null by Lemma~\ref{lem pn-free}. Thus there exists a height one prime ideal $\pe\subseteq \Lambda$ such that $N_\pe$ is not trivial. Since $\chara_\Lambda X_\Sigma(\wtk)$ is a prime ideal, it follows that $\pe=\chara_\Lambda X_\Sigma(\wtk)$ and $\length_{\Lambda_\pe}X_\Sigma(\wtk)_\pe=1$. Hence $N_\pe=X_\Sigma(\wtk)_\pe$, so $X(\wtk)_\pe$ is trivial. On the other hand, for a height one prime ideal $\qu\subseteq \Lambda$ different from $\chara_\Lambda X_\Sigma(\wtk)$, we have $X(\wtk)_\qu=0$ because $X_\Sigma(\wtk)_\qu=0$. Therefore, we conclude that $X(\wtk)$ is pseudo-null. 
\end{proof}

\subsection{Proof of Theorem~\ref{thm NC wGC to NF}}
\begin{proof}[Proof of Theorem~\ref{thm NC wGC to NF}]
To show the contraposition, assume that wGC is true for $(k^+, p)$, and $\calG^{(2)}(k^\cy_\infty)$ is a free metabelian pro-$p$-group. Then the module $X(k^\cy_\infty)$, which is the abelianization of $\calG^{(2)}(k^\cy_\infty)$, is torsion-free as a $\Z_p$-module. Thus we see that $X(k^+_\infty)\simeq X(k^\cy_\infty)^+$ is trivial, since we assume the validity of wGC for $(k^+,p)$. Hence $L(k^\cy_\infty)=M_\Sigma(k^\cy_\infty)$ holds by Lemma~\ref{lem LM}. 

Since $\widetilde{k}/k^\cy_\infty$ is unramified, it follows that $L(\wtk)/k^\cy_\infty$ is an unramified metabelian pro-$p$-extension. Thus there is a natural surjection $\calG^{(2)}(k^\cy_\infty)\to \Gal(L(\wtk)/k^\cy_\infty)$. On the other hand, there is a natural surjection $\Gal(M_\Sigma(\widetilde{k})/k^\cy_\infty)\to\Gal(L(\wtk)/k^\cy_\infty)$. Since $\Gal(M_\Sigma(\wtk)/k^\cy_\infty)$ and $\Gal(L(\wtk)/k^\cy_\infty)$ are metabelian pro-$p$-groups, by the freeness of $\calG^{(2)}(k^\cy_\infty)$, there exists a morphism $\varphi\colon\calG^{(2)}(k^\cy_\infty)\to \Gal(M_\Sigma(\wtk)/k^\cy_\infty)$ that satisfies the following commutative diagram. 
\[
\begin{tikzpicture}
\node (a) at (0,2) {$\calG^{(2)}(k^\cy_\infty)$}; 
\node (b) at (0,0) {$\Gal(M_\Sigma(\wtk)/k^\cy_\infty)$}; 
\node (c) at (4.3,0) {$\Gal(L(\wtk)/k^\cy_\infty)$}; 

\draw[->] (a) -- node[auto][swap]{$\scriptstyle \varphi$} (b); \draw[->] (b) -- (c); \draw[->] (a) -- (c); 
\end{tikzpicture}
\]
By taking the abelianizations of the diagram, we obtain the following commutative diagram. 
\[
\begin{tikzpicture}
\node (a) at (0,2) {$X(k^\cy_\infty)$}; 
\node (b) at (0,0) {$X_\Sigma(k^\cy_\infty)$}; 
\node (c) at (4,0) {$X(k^\cy_\infty)$}; 

\draw[->] (a) -- node[auto][swap]{$\scriptstyle \varphi_{\mathrm{ab}}$} (b); \draw[->] (b) -- node[auto]{$\scriptstyle =$} (c); \draw[->] (a) -- node[auto]{$\scriptstyle =$} (c); 
\end{tikzpicture}
\]
Here, the bottom horizontal map is the identity by Lemma~\ref{lem LM}. We see that $\varphi_{\mathrm{ab}}$ is the identity and, in particular, a surjection, so $\varphi$ is surjective. By the commutativity of the first diagram, we see that $\varphi$ maps $\Gal(\calL^{(2)}(k^\cy_\infty)/\wtk)$ into $X_\Sigma(\wtk)$, which is abelian. Thus $D_1(\Gal(\calL^{(2)}(k^\cy_\infty)/\wtk))=\Gal(\calL^{(2)}(k^\cy_\infty)/L(\wtk))$ is contained in the kernel of $\varphi$. Hence it follows that $\varphi$ induces a surjective morphism $\overline{\varphi}\colon\Gal(L(\wtk)/k^\cy_\infty)\to \Gal(M_\Sigma(\wtk)/k^\cy_\infty)$ that satisfies the following commutative diagram. 
\[
\begin{tikzpicture}
\node (a) at (0,2) {$\Gal(L(\wtk)/k^\cy_\infty)$}; 
\node (b) at (0,0) {$\Gal(M_\Sigma(\wtk)/k^\cy_\infty)$}; 
\node (c) at (4.3,0) {$\Gal(L(\wtk)/k^\cy_\infty)$}; 

\draw[->] (a) -- node[auto][swap]{$\scriptstyle \overline{\varphi}$} (b); \draw[->] (b) -- (c); \draw[->] (a) -- node[auto]{$\scriptstyle =$} (c); 
\end{tikzpicture}
\]
The morphism $\overline{\varphi}$ is an isomorphism by the commutativity, so the natural surjection $\Gal(M_\Sigma(\wtk)/k^\cy_\infty)\to \Gal(L(\wtk)/k^\cy_\infty)$ is an isomorphism, as well. This means that $L(\wtk)=M_\Sigma(\wtk)$ holds, as desired. 
\end{proof}

\section{Proof of Theorem~\ref{thm main}}\label{section proof main}
Let $k$ be a quartic CM-field, $p$ an odd prime number that splits completely in $k$, and $\Sigma$ a $p$-adic CM-type of $k$. Suppose that $X(k^\cy_\infty)\simeq \Z_p^2$ as a $\Z_p\bbr{\Gal(k^\cy_\infty/k)}$-module. We keep this setting in the present section. 

The implication \ref{main item 1}$\Rightarrow$\ref{main item 2} is clear. The implication \ref{main item 2}$\Rightarrow$\ref{main item 3} follows from Proposition~\ref{prop wGGC different}. The implication \ref{main item 3}$\Rightarrow$\ref{main item 4} follows from Theorem~\ref{thm NC wGC to NF} combined with the following lemma. 
\begin{lem}\label{lem trivial}
$X(k^+_\infty)$ is trivial. 
\end{lem}
\begin{proof}
It is mentioned in \cite[Remark~4.1]{Yanagisawa2026}. Note that $p$-adic Leopoldt's Conjecture is true for $k$, because the $\Z$-rank of the unit group of $k$ is $1$. By \cite[Proposition~1]{Ozaki1997}, there is a natural surjection $X(k^\cy_\infty)\to \Gal(\wtk/k^\cy_\infty)$ of $\Z_p\bbr{\Gal(k^\cy_\infty/k)}$-modules. By the assumption $X(k^\cy_\infty)\simeq \Z_p^2$, we see that the surjection is an isomorphism. Hence $X(k^+_\infty)\simeq X(k^\cy_\infty)^+=0$ by the validity of $p$-adic Leopoldt's Conjecture. 
\end{proof}
\subsection{Proof of the implication \ref{main item 3}$\Rightarrow$\ref{main item 1}}
We mentioned the cyclicity of $X_\Sigma(\wtk)$ as a $\Z_p\bbr{\Gal(\wtk/k^\cy_\infty)}$-module in \cite[Remark~5.4]{Yanagisawa2026}. Here we show a refined statement. 
\begin{lem}\label{lem free rank 1}
$X_\Sigma(\wtk)$ is a free $\Z_p\bbr{\Gal(\wtk/k^\cy_\infty)}$-module of rank $1$. 
\end{lem}
\begin{proof}
Since $X(k^+_\infty)$ is trivial by Lemma~\ref{lem trivial}, it follows from Lemma~\ref{lem four-term} that $X_\Sigma(\wtk)_{\Gal(\wtk/k^\cy_\infty)}$ is cyclic as a $\Z_p$-module. Hence $X_\Sigma(\wtk)$ is cyclic as a $\Z_p\bbr{\Gal(\wtk/k^\cy_\infty)}$-module by Nakayama's lemma. Since $X_\Sigma(\wtk)$ is not pseudo-null as a $\Lambda$-module by Proposition~\ref{prop non-trivial} and Lemma~\ref{lem pn-free}, we see that $X_\Sigma(\wtk)$ is torsion-free as a $\Z_p\bbr{\Gal(\wtk/k^\cy_\infty)}$-module by Lemma~\ref{lem tor pn}. Hence $X_\Sigma(\wtk)$ is free of rank $1$ as a $\Z_p\bbr{\Gal(\wtk/k^\cy_\infty)}$-module.  
\end{proof}
\begin{lem}\label{lem prime element}
There exists a prime element $\pi\in\Lambda$ such that $X_\Sigma(\wtk)$ is isomorphic to $\Lambda/(\pi)$ as a $\Lambda$-module. 
\end{lem}
\begin{proof}
We can take a basis $x$ of $X_\Sigma(\wtk)$ as a $\Z_p\bbr{\Gal(\wtk/k^\cy_\infty)}$-module by Lemma~\ref{lem free rank 1}. For each $f\in\Lambda$, there exists a unique element $\varphi(f)\in\Z_p\bbr{\Gal(\wtk/k^\cy_\infty)}$ such that $f x=\varphi(f) x$. This correspondence defines a surjective ring homomorphism $\varphi\colon\Lambda\to \Z_p\bbr{\Gal(\wtk/k^\cy_\infty)}$. Then, if we let $\pe$ be the kernel of $\varphi$, we see that $\pe$ is a height one prime ideal of $\Lambda$, so there is a prime element $\pi$ that generates $\pe$. Here we note that the correspondence $f\mapsto \varphi(f)x$ defines a surjective morphism $\Lambda\to X_\Sigma(\wtk)$ of $\Lambda$-modules. This shows $\Lambda/(\pi)\simeq X_\Sigma(\wtk)$, as desired. 
\end{proof}
By Lemma~\ref{lem prime element}, we see that $\chara_\Lambda X_\Sigma(\wtk)$ is a prime ideal. Thus the implication \ref{main item 3}$\Rightarrow$\ref{main item 1} follows from Theorem~\ref{thm general}.

\subsection{Proof of the implication \ref{main item 4}$\Rightarrow$\ref{main item 3}}
\begin{lem}\label{lem D1}
There is a natural isomorphism $D_1(\calG(k^\cy_\infty))_\ab\simeq X(\wtk)$ of $\Z_p\bbr{\Gal(\wtk/k^\cy_\infty)}$-modules. 
\end{lem}
\begin{proof}
Since $L(k^\cy_\infty)=\wtk$ holds by the assumption, we have $D_1(\calG(k^\cy_\infty))=\Gal(\calL(k^\cy_\infty)/\wtk)$. On the other hand, since $\wtk/k^\cy_\infty$ is unramified, we have $\calL(\wtk)=\calL(k^\cy_\infty)$, so $D_1(\calG(k^\cy_\infty))=\calG(\wtk)$ holds. Then take the abelianization. 
\end{proof}

Take a minimal presentation
\[
\begin{tikzpicture}
\node (a) at (0,0) {$1$}; 
\node (b) at (1.5,0) {$R$}; 
\node (c) at (3,0) {$F$}; 
\node (d) at (4.8,0) {$\calG(k^\cy_\infty)$}; 
\node (e) at (6.6,0) {$1$}; 

\draw[->] (a) -- (b); \draw[->] (b) -- (c); \draw[->] (c) -- (d); \draw[->] (d) -- (e); 
\end{tikzpicture}
\]
of $\calG(k^\cy_\infty)$ as a pro-$p$-group. Here we note that $F$ is a free pro-$p$-group of rank $2$, because $\calG(k^\cy_\infty)_\ab\simeq \Z_p^2$ holds by the assumption. So the natural surjective morphism $F_\ab\to\calG(k^\cy_\infty)_\ab$ is an isomorphism, because both are isomorphic to $\Z_p^2$. 
\begin{lem}\label{lem kernel}
The kernel of the natural morphism $D_1(F)_\ab\to D_1(\calG(k^\cy_\infty))_\ab$ is $RD_2(F)/D_2(F)$. 
\end{lem}
\begin{proof}
Since the natural surjection $F_\ab\to \calG(k^\cy_\infty)_\ab$ is an isomorphism as we mentioned above, it follows that $R\subseteq D_1(F)$ holds. Thus there is an exact sequence
\[
\begin{tikzpicture}
\node (a0) at (0,0) {$1$}; 
\node (b0) at (1.7,0) {$R$}; 
\node (c0) at (3.9,0) {$D_1(F)$}; 
\node (d0) at (6.8,0) {$D_1(\calG(k^\cy_\infty))$}; 
\node (e0) at (9,0) {$1$}; 

\draw[->] (a0) -- (b0); \draw[->] (b0) -- (c0); \draw[->] (c0) -- (d0); \draw[->] (d0) -- (e0); 
\end{tikzpicture}
\]
of pro-$p$-groups. Then the assertion follows from the right exactness of abelianization. 
\end{proof}
\begin{proof}[Proof of the implication \ref{main item 4}$\Rightarrow$\ref{main item 3}]
To show the contraposition, suppose that $L(\wtk)=M_\Sigma(\wtk)$ holds. Then there is a natural isomorphism $D_1(\calG(k^\cy_\infty))_\ab\simeq X_\Sigma(\wtk)$ of $\Z_p\bbr{\Gal(\wtk/k^\cy_\infty)}$-modules by Lemma~\ref{lem D1}. Moreover, it follows that $D_1(\calG(k^\cy_\infty))_\ab$ is free of rank $1$ as a $\Z_p\bbr{\Gal(\wtk/k^\cy_\infty)}$-module by Lemma~\ref{lem free rank 1}. 
On the other hand, since there is a natural isomorphism $F_\ab\simeq\Gal(\wtk/k^\cy_\infty)$, it can also be shown that $D_1(F)_\ab$ is a free $\Z_p\bbr{\Gal(\wtk/k^\cy_\infty)}$-module of rank $1$ by Lemma~\ref{lem rank 2 free}. Thus the surjection $D_1(F)_\ab\to D_1(\calG(k^\cy_\infty))_\ab$ is an isomorphism of $\Z_p\bbr{\Gal(\wtk/k^\cy_\infty)}$-modules. Hence we obtain $R\subseteq D_2(F)$ by Lemma~\ref{lem kernel}, so $F/D_2(F)\simeq \calG(k^\cy_\infty)/D_2(\calG(k^\cy_\infty))\simeq \calG^{(2)}(k^\cy_\infty)$ holds. Therefore, we conclude that $\calG^{(2)}(k^\cy_\infty)$ is a free metabelian pro-$p$-group. 
\end{proof}
This completes the proof of Theorem~\ref{thm main}.

\section{A sufficient condition for Greenberg's Generalized Conjecture}\label{section criterion}
By using Theorem~\ref{thm main}, we give a new sufficient condition for GGC over a quartic CM-field. Our strategy is to give a sufficient condition for Condition~\ref{main item 4} in Theorem~\ref{thm main}. 

For a finite extension $K/\Q$, we let $h(K)$ denote the class number of $K$ and $v_p(h(K))$ its $p$-order. Let $K^\cy_1/K$ denote the unique subextension of $K^\cy_\infty/K$ of degree $p$. 
\begin{prop}\label{prop criterion}
Let $k$ be a quartic CM-field and $p$ an odd prime number that splits completely in $k$. Suppose that $X(k^\cy_\infty)\simeq \Z_p^2$ as a $\Z_p\bbr{\Gal(k^\cy_\infty/k)}$-module. If there exists an unramified cyclic extension $K/k$ of degree $p$ that satisfies $v_p(h(K^\cy_1))-v_p(h(K))\leq p$, then $\calG^{(2)}(k^\cy_\infty)$ is not a free metabelian pro-$p$-group, so GGC is true for $(k, p)$ by Theorem~\ref{thm main}. 
\end{prop}
\begin{proof}
To show the contraposition, suppose that $\calG^{(2)}(k^\cy_\infty)$ is a free metabelian pro-$p$-group. Then it suffices to show that an inequality $v_p(h(K^\cy_1))-v_p(h(K))\geq p+1$ holds for arbitrary unramified cyclic extension $K/k$ of degree $p$. 

Since $L(k^\cy_\infty)=\wtk$ holds by the assumption, we have $L(k)\subseteq \wtk$, so $K\subseteq \wtk$ holds. Then $K^\cy_\infty/k^\cy_\infty$ is a subextension of $\wtk/k^\cy_\infty$ of degree $p$. Note that $\widetilde{k}/k^\cy_\infty$ is unramified by \cite[Proposition~1]{Ozaki1997}, hence $\wtk/K^\cy_\infty$ is unramified as well. Since $K/\Q$ is unramified at $p$, it follows that $K^\cy_\infty/K$ is totally ramified at every $p$-adic prime. 

We show $X(K^\cy_\infty)\simeq \Z_p^{p+1}$, in what follows. Take a minimal presentation
\[
\begin{tikzpicture}
\node (a) at (0,0) {$1$}; 
\node (b) at (1.5,0) {$R$}; 
\node (c) at (3,0) {$F$}; 
\node (d) at (4.8,0) {$\calG(k^\cy_\infty)$}; 
\node (e) at (6.6,0) {$1$}; 

\draw[->] (a) -- (b); \draw[->] (b) -- (c); \draw[->] (c) -- (d); \draw[->] (d) -- (e); 
\end{tikzpicture}
\]
of $\calG(k^\cy_\infty)$ as a pro-$p$-group. Since $\calG^{(2)}(k^\cy_\infty)$ is a free metabelian pro-$p$-group, we have $R\subseteq D_2(F)$. Let $F_0$ be the inverse image of $\calG(K^\cy_\infty)\subseteq \calG(k^\cy_\infty)$ in $F$. Then $F_0$ is an open subgroup of $F$ of index $p$, so it is a free pro-$p$-group of rank $p+1$ by Nielsen--Schreier theorem (see \cite[Theorem~3.6.2~(b)]{RZ2010}). Since $D_1(F)\subseteq F_0$ holds, we have $D_2(F)\subseteq D_1(F_0)$, so $R\subseteq D_1(F_0)$ holds. Thus there is an exact sequence
\[
\begin{tikzpicture}
\node (a) at (0,0) {$1$}; 
\node (b) at (1.5,0) {$R$}; 
\node (c) at (3,0) {$F_0$}; 
\node (d) at (4.8,0) {$\calG(K^\cy_\infty)$}; 
\node (e) at (6.6,0) {$1$}; 

\draw[->] (a) -- (b); \draw[->] (b) -- (c); \draw[->] (c) -- (d); \draw[->] (d) -- (e); 
\end{tikzpicture}
\]
of pro-$p$-groups. By taking the abelianizations, we obtain $X(K^\cy_\infty)\simeq (F_0)_\ab\simeq \Z_p^{p+1}$, as desired. 

Take a topological generator $\gamma\in\Gal(K^\cy_\infty/K)$ and define $\nu\coloneqq 1+\gamma+\cdots+\gamma^{p-1}$. Put $Y\coloneqq\Ker(X(K^\cy_\infty)\to X(K))=\Gal(L(K^\cy_\infty)/K^\cy_\infty L(K))$. Then, by the proof of Iwasawa's class number formula (\cite[\S 13.3]{Washington1997}), there are natural isomorphisms $X(K^\cy_\infty)/Y\simeq X(K)$ and $X(K^\cy_\infty)/\nu Y\simeq X(K^\cy_1)$, so it suffices to show an inequality $\#(Y/\nu Y)\geq p^{p+1}$. Put $Z\coloneqq\Gal(L(K^\cy_\infty)/\wtk L(K))$, so that we have $Y/Z\simeq \Gal(\wtk/\wtk\cap K^\cy_\infty L(K))$. Since $\wtk/K$ is abelian, we see that $\gamma$ acts trivially on $\Gal(\wtk/\wtk\cap K^\cy_\infty L(K))$, hence $\nu$ acts on $\Gal(\wtk/\wtk\cap K^\cy_\infty L(K))$ as multiplication by $p$. Thus there is an exact commutative diagram
\[
\begin{tikzpicture}
\node (a) at (0,2) {$0$}; 
\node (b) at (1.7,2) {$Z$}; 
\node (c) at (3.5,2) {$Y$}; 
\node (d) at (6.7,2) {$\Gal(\wtk/\wtk\cap K^\cy_\infty L(K))$}; 
\node (e) at (9.5,2) {$0$}; 

\draw[->] (a) -- (b); \draw[->] (b) -- (c); \draw[->] (c) -- (d); \draw[->] (d) -- (e); 

\node (a0) at (0,0) {$0$}; 
\node (b0) at (1.7,0) {$Z$}; 
\node (c0) at (3.5,0) {$Y$}; 
\node (d0) at (6.7,0) {$\Gal(\wtk/\wtk\cap K^\cy_\infty L(K))$}; 
\node (e0) at (9.5,0) {$0$}; 

\draw[->] (a0) -- (b0); \draw[->] (b0) -- (c0); \draw[->] (c0) -- (d0); \draw[->] (d0) -- (e0); 

\draw[->] (b) -- node[auto][swap]{$\scriptstyle \times\nu$} (b0); \draw[->] (c) -- node[auto][swap]{$\scriptstyle \times\nu$} (c0); \draw[->] (d) -- node[auto][swap]{$\scriptstyle \times p$} (d0); 
\end{tikzpicture}
\]
of $\Z_p\bbr{\Gal(K^\cy_\infty/K)}$-modules. By Snake lemma, we obtain $\#(Y/\nu Y)=p^2 \#(Z/\nu Z)$. Since $Z$ is isomorphic to $\Z_p^{p-1}$, we have $\nu (Z/pZ)=0$. Indeed, since $Z/pZ$ is a $\F_p$-vector space of dimension $p-1$, there is an integer $i$ with $1\leq i\leq p-1$ and a decreasing filtration
\[
Z/pZ\supsetneq (\gamma-1)(Z/pZ) \supsetneq (\gamma-1)^2(Z/pZ)\supsetneq\cdots\supsetneq (\gamma-1)^i(Z/pZ)=0, 
\]
so we have $(\gamma-1)^{p-1}(Z/pZ)=0$. Here we have 
\[
\nu=\frac{\gamma^p-1}{\gamma-1}\equiv \frac{(\gamma-1)^p}{\gamma-1}=(\gamma-1)^{p-1} \pmod{p}, 
\]
hence we obtain $\nu (Z/pZ)=0$, as desired. Thus there is a natural surjection $Z/\nu Z\to Z/pZ\simeq \F_p^{p-1}$. Therefore, we obtain $\#(Y/\nu Y)\geq p^2 p^{p-1}=p^{p+1}$, as desired. 
\end{proof}

To provide non-trivial examples, we use the following lemma. 
\begin{lem}[{\cite[Lemma~4.10]{Yanagisawa2026}}]\label{lem zp2}
Let $k$ be a quartic CM-field and $p$ an odd prime number that splits completely in $k$. Suppose that $\rank_p X(k)=\rank_p X(k^\cy_1)=2$ holds. Then $X(k^\cy_\infty)\simeq \Z_p^2$ as a $\Z_p\bbr{\Gal(k^\cy_\infty/k)}$-module and $X(\wtk)$ is not trivial. 
\end{lem}
By using Proposition~\ref{prop criterion} and Lemma~\ref{lem zp2}, we provide examples of GGC. The following examples are checked by using PARI/GP, under the assumption of Generalized Riemann Hypothesis.

\begin{example}
Let $p=3$ and $k=\Q(\sqrt{-2834-654\sqrt{13}})$. Then $k$ is an imaginary abelian quartic field such that $\Gal(k/\Q)\simeq \Z/4\Z$. Let $K$ be the field defined by $x^{12} - 2x^{11} - 52x^{10} + 237x^9 + 654x^8 - 6792x^7 + 18999x^6 - 25662x^5 + 40506x^4 - 106393x^3 + 172630x^2 - 134848x + 41011$ over $\Q$. Then $K/k$ is an unramified cyclic extension of degree $3$. It can be checked that the assumptions of Lemma~\ref{lem zp2} and Proposition~\ref{prop criterion} are satisfied. Hence $X(\wtk)$ is not trivial and GGC is true for $(k, 3)$.  
\end{example}

\begin{example}
Let $p=3$ and $k=\Q(\sqrt{-1514+234\sqrt{13}})$. Then $k$ is a quartic CM-field that is not a Galois extension over $\Q$. Let $K$ be the field defined by $x^{12} + 10x^{10} + 31x^8 + 787x^6 + 5315x^4 + 10994x^2 + 98773$ over $\Q$. Then $K/k$ is an unramified cyclic extension of degree $3$. It can be checked that the assumptions of Lemma~\ref{lem zp2} and Proposition~\ref{prop criterion} are satisfied. Hence $X(\wtk)$ is not trivial and GGC is true for $(k, 3)$.  
\end{example}
As we mentioned in Remark~\ref{rem Bandini}, our method does not give new examples of GGC, in the case where $k$ is an imaginary biquadratic field. 

In the setting of Proposition~\ref{prop criterion}, if the class number of $k$ is not divisible by $p$, then GGC follows from Fujii's result \cite[Theorem~2]{Fujii2017} (see \cite[Theorem~1.1]{Itoh2011} as well). We note that the above examples are not covered by Fujii's result.

\subsection*{Use of AI}
The author used ChatGPT (GPT-6, OpenAI) for the mathematical discussions during the preparation of the present article. It was used to obtain proof ideas of the main results. ChatGPT was also used to check the manuscript for grammatical and typographical errors. 

All mathematical arguments are independently verified by the author and all parts of the manuscript were written by the author.

\subsection*{Acknowledgement}
I would like to thank my supervisor Takenori Kataoka for reading the manuscript and providing valuable comments.

\printbibliography

@article {Iwasawa1973,
    AUTHOR = {Iwasawa, Kenkichi},
     TITLE = {On {${\bf Z}\sb{l}$}-extensions of algebraic number fields},
   JOURNAL = {Ann. of Math. (2)},
  FJOURNAL = {Annals of Mathematics. Second Series},
    VOLUME = {98},
      YEAR = {1973},
     PAGES = {246--326},
}

@article {Greenberg1976,
    AUTHOR = {Greenberg, Ralph},
     TITLE = {On the {I}wasawa invariants of totally real number fields},
   JOURNAL = {Amer. J. Math.},
  FJOURNAL = {American Journal of Mathematics},
    VOLUME = {98},
      YEAR = {1976},
    NUMBER = {1},
     PAGES = {263--284},
}

@article {Greenberg1978,
    AUTHOR = {Greenberg, Ralph},
     TITLE = {On the structure of certain {G}alois groups},
   JOURNAL = {Invent. Math.},
  FJOURNAL = {Inventiones Mathematicae},
    VOLUME = {47},
      YEAR = {1978},
    NUMBER = {1},
     PAGES = {85--99},
}

@article {Wintenberger1980,
    AUTHOR = {Wintenberger, Jean-Pierre},
     TITLE = {Structure galoisienne de limites projectives d'unit\'es
              locales},
   JOURNAL = {Compositio Math.},
  FJOURNAL = {Compositio Mathematica},
    VOLUME = {42},
      YEAR = {1980/81},
    NUMBER = {1},
     PAGES = {89--103},
}

@incollection {Nguyen1984,
    AUTHOR = {Nguyen Quang Do, Thong},
     TITLE = {Formations de classes et modules d'{I}wasawa},
 BOOKTITLE = {Number theory, {N}oordwijkerhout 1983 ({N}oordwijkerhout,
              1983)},
    SERIES = {Lecture Notes in Math.},
    VOLUME = {1068},
     PAGES = {167--185},
 PUBLISHER = {Springer, Berlin},
      YEAR = {1984},
}

@article {HT1994,
    AUTHOR = {Hida, H. and Tilouine, J.},
     TITLE = {On the anticyclotomic main conjecture for {CM} fields},
   JOURNAL = {Invent. Math.},
  FJOURNAL = {Inventiones Mathematicae},
    VOLUME = {117},
      YEAR = {1994},
    NUMBER = {1},
     PAGES = {89--147},
}

@book {Washington1997,
    AUTHOR = {Washington, Lawrence C.},
     TITLE = {Introduction to cyclotomic fields},
    SERIES = {Graduate Texts in Mathematics},
    VOLUME = {83},
   EDITION = {Second},
 PUBLISHER = {Springer-Verlag, New York},
      YEAR = {1997},
}

@article {Ozaki1997,
    AUTHOR = {Ozaki, Manabu},
     TITLE = {The class group of {${\bf Z}_p$}-extensions over totally real
              number fields},
   JOURNAL = {Tohoku Math. J. (2)},
  FJOURNAL = {The Tohoku Mathematical Journal. Second Series},
    VOLUME = {49},
      YEAR = {1997},
    NUMBER = {3},
     PAGES = {431--435},
}

@incollection {Greenberg2001,
    AUTHOR = {Greenberg, Ralph},
     TITLE = {Iwasawa theory---past and present},
 BOOKTITLE = {Class field theory---its centenary and prospect ({T}okyo,
              1998)},
    SERIES = {Adv. Stud. Pure Math.},
    VOLUME = {30},
     PAGES = {335--385},
 PUBLISHER = {Math. Soc. Japan, Tokyo},
      YEAR = {2001},
}

@article {Nguyen2006,
    AUTHOR = {Nguyen Quang Do, Thong},
     TITLE = {Sur la conjecture faible de {G}reenberg dans le cas ab\'elien
              {$p$}-d\'ecompos\'e},
   JOURNAL = {Int. J. Number Theory},
  FJOURNAL = {International Journal of Number Theory},
    VOLUME = {2},
      YEAR = {2006},
    NUMBER = {1},
     PAGES = {49--64},
}

@article {Ozaki2007,
    AUTHOR = {Ozaki, Manabu},
     TITLE = {Non-abelian {I}wasawa theory of {$\Bbb Z_p$}-extensions},
   JOURNAL = {J. Reine Angew. Math.},
  FJOURNAL = {Journal f\"ur die Reine und Angewandte Mathematik. [Crelle's
              Journal]},
    VOLUME = {602},
      YEAR = {2007},
     PAGES = {59--94},
}

@article {Bandini2007,
    AUTHOR = {Bandini, Andrea},
     TITLE = {Greenberg's conjecture and capitulation in {$\Bbb
              Z_p^d$}-extensions},
   JOURNAL = {J. Number Theory},
  FJOURNAL = {Journal of Number Theory},
    VOLUME = {122},
      YEAR = {2007},
    NUMBER = {1},
     PAGES = {121--134},
}

@book {NSW2008,
    AUTHOR = {Neukirch, J\"urgen and Schmidt, Alexander and Wingberg, Kay},
     TITLE = {Cohomology of number fields},
    SERIES = {Grundlehren der mathematischen Wissenschaften [Fundamental
              Principles of Mathematical Sciences]},
    VOLUME = {323},
   EDITION = {Second},
 PUBLISHER = {Springer-Verlag, Berlin},
      YEAR = {2008},
     PAGES = {xvi+825},
}

@book {RZ2010,
    AUTHOR = {Ribes, Luis and Zalesskii, Pavel},
     TITLE = {Profinite groups},
    SERIES = {Ergebnisse der Mathematik und ihrer Grenzgebiete. 3. Folge. A
              Series of Modern Surveys in Mathematics [Results in
              Mathematics and Related Areas. 3rd Series. A Series of Modern
              Surveys in Mathematics]},
    VOLUME = {40},
   EDITION = {Second},
 PUBLISHER = {Springer-Verlag, Berlin},
      YEAR = {2010},
}

@article {HO2010,
    AUTHOR = {Hachimori, Yoshitaka and Ochiai, Tadashi},
     TITLE = {Notes on non-commutative {I}wasawa theory},
   JOURNAL = {Asian J. Math.},
  FJOURNAL = {Asian Journal of Mathematics},
    VOLUME = {14},
      YEAR = {2010},
    NUMBER = {1},
     PAGES = {11--17},
}

@article {Itoh2011,
    AUTHOR = {Itoh, Tsuyoshi},
     TITLE = {On multiple {$\Bbb Z_p$}-extensions of imaginary abelian
              quartic fields},
   JOURNAL = {J. Number Theory},
  FJOURNAL = {Journal of Number Theory},
    VOLUME = {131},
      YEAR = {2011},
    NUMBER = {1},
     PAGES = {59--66},
}

@article {Fujii2011,
    AUTHOR = {Fujii, Satoshi},
     TITLE = {On the depth of the relations of the maximal unramified
              pro-{$p$} {G}alois group over the cyclotomic {$\Bbb
              Z_p$}-extension},
   JOURNAL = {Acta Arith.},
  FJOURNAL = {Acta Arithmetica},
    VOLUME = {149},
      YEAR = {2011},
    NUMBER = {2},
     PAGES = {101--110},
}

@incollection {Fujii2011RIMS,
    AUTHOR = {Fujii, Satoshi},
     TITLE = {Greenberg's generalized conjecture and unramified {G}alois
              groups over the cyclotomic {$\Bbb Z_p$}-extensions},
 BOOKTITLE = {Algebraic number theory and related topics 2009},
    SERIES = {RIMS K\^oky\^uroku Bessatsu},
    VOLUME = {B25},
     PAGES = {95--101},
 PUBLISHER = {Res. Inst. Math. Sci. (RIMS), Kyoto},
      YEAR = {2011},
}

@article {Fujii2014,
    AUTHOR = {Fujii, Satoshi},
     TITLE = {On restricted ramifications and pseudo-null submodules of
              {I}wasawa modules for {$\Bbb Z_p^2$}-extensions},
   JOURNAL = {J. Ramanujan Math. Soc.},
  FJOURNAL = {Journal of the Ramanujan Mathematical Society},
    VOLUME = {29},
      YEAR = {2014},
    NUMBER = {3},
     PAGES = {295--305},
}

@article {Fujii2017,
    AUTHOR = {Fujii, Satoshi},
     TITLE = {On {G}reenberg's generalized conjecture for {CM}-fields},
   JOURNAL = {J. Reine Angew. Math.},
  FJOURNAL = {Journal f\"ur die Reine und Angewandte Mathematik. [Crelle's
              Journal]},
    VOLUME = {731},
      YEAR = {2017},
     PAGES = {259--278},
}

@article {Nguyen2017,
    AUTHOR = {Nguyen Quang Do, Thong},
     TITLE = {Sur une forme faible de la conjecture de {G}reenberg {II}},
   JOURNAL = {Int. J. Number Theory},
  FJOURNAL = {International Journal of Number Theory},
    VOLUME = {13},
      YEAR = {2017},
    NUMBER = {4},
     PAGES = {1061--1070},
}

@article {Fujii2022,
    AUTHOR = {Fujii, Satoshi},
     TITLE = {On families of imaginary abelian fields with pseudo-null
              unramified {I}wasawa modules},
   JOURNAL = {New York J. Math.},
  FJOURNAL = {New York Journal of Mathematics},
    VOLUME = {28},
      YEAR = {2022},
     PAGES = {523--533},
}

@article {Murakami2023,
    AUTHOR = {Murakami, Kazuaki},
     TITLE = {A weak form of {G}reenberg's generalized conjecture for
              imaginary quadratic fields},
   JOURNAL = {J. Number Theory},
  FJOURNAL = {Journal of Number Theory},
    VOLUME = {244},
      YEAR = {2023},
     PAGES = {308--338},
}

@article {Jaulent2024,
    AUTHOR = {Jaulent, Jean-Fran\c cois},
     TITLE = {Sur la trivialit\'e{} de certains modules d'{I}wasawa},
   JOURNAL = {Funct. Approx. Comment. Math.},
  FJOURNAL = {Uniwersytet im. Adama Mickiewicza w Poznaniu. Wydzia\l\
              Matematyki i Informatyki. Functiones et Approximatio
              Commentarii Mathematici},
    VOLUME = {70},
      YEAR = {2024},
    NUMBER = {1},
     PAGES = {29--39},
}

@misc{Okano2026pre,
      title={On the triviality of the unramified Iwasawa modules of the maximal multiple $\mathbb{Z}_p$-extensions}, 
      author={Keiji Okano},
      year={2026},
      eprint={2406.01892v5},
      archivePrefix={arXiv},
}

@misc{MO2026pre,
      title={Trivial Iwasawa modules and arithmetic linking numbers}, 
      author={Yasushi Mizusawa and Keiji Okano},
      year={2026},
      eprint={2610.31965v1},
      archivePrefix={arXiv},
}

@article {Yanagisawa2026,
    AUTHOR = {Yanagisawa, Takuya},
     TITLE = {On a bound of {$p$}-ranks of {I}wasawa modules of
              {$\Bbb{Z}_{p}$}-extensions over a quartic {CM}-field},
   JOURNAL = {J. Number Theory},
  FJOURNAL = {Journal of Number Theory},
    VOLUME = {289},
      YEAR = {2026},
     PAGES = {96--113},
}

@unpublished{WingbergPre,
  author = {Wingberg, Kay},
  title  = {Free pro-$p$ extensions of number fields},
  note   = {preprint},
}
\end{document}